\documentclass[letterpaper]{article}

\usepackage{amsmath,amsfonts,amsthm}
\usepackage{graphicx}

\newtheorem{defs}{Definition}
\newtheorem{lem}{Lemma}

\newtheorem{conj}{Conjecture}

\title{Proof of the volume conjecture for Whitehead chains}
\author{Roland van der Veen}
\begin{document}

\maketitle

\begin{abstract}
\noindent We prove the volume conjecture for an infinite family of links called
Whitehead chains that generalizes both the Whitehead link and the
Borromean rings.

\end{abstract}

\section{Introduction}

The volume conjecture relates the colored Jones polynomials of a
knot to the simplicial volume of its complement. More precisely,
let us denote the normalized $N$-colored Jones polynomial of a
knot $K$ by $J_N(K)$ and let $\mathrm{Vol}(K)$ be $v_3$ times the
simplicial volume (Gromov norm) of $\mathbb{S}^3 - K$, where $v_3$
is the volume of the hyperbolic regular ideal tetrahedron. The
volume conjecture can now be stated as follows:

\begin{conj} (\textbf{Volume conjecture}) \cite{Murakami2}\\
For any knot $K$ we have: $$\lim_{N\to\infty}\frac{2\pi}{N}\log|J_N(K)(e^{\frac{2\pi i}{N}})| = \mathrm{Vol}(K)$$
\end{conj}

\noindent So far the conjecture has been proven for only the
figure eight knot, torus knots, Whitehead doubles of certain torus
knots and connected sums of these knots, \cite{Murakamiintro},
\cite{KashaevTirkkonen}, \cite{HaoZheng}.

It is well known that the volume conjecture is false for many
splittable links so it is unclear how to extend the volume
conjecture to links. On the other hand, the volume conjecture has
been shown to hold for the Whitehead link \cite{HaoZheng} and the
Borromean rings \cite{Garoufalidis}. In this paper, we introduce
the family of Whitehead chains generalizing both the Whitehead
link and the Borromean rings and we settle the volume conjecture
for this family.

The Whitehead chains are defined in terms of the tangles Belt,
Clasp and Twist depicted in figure 1 on the next page.

\newpage

\begin{figure}
\begin{center}
\includegraphics{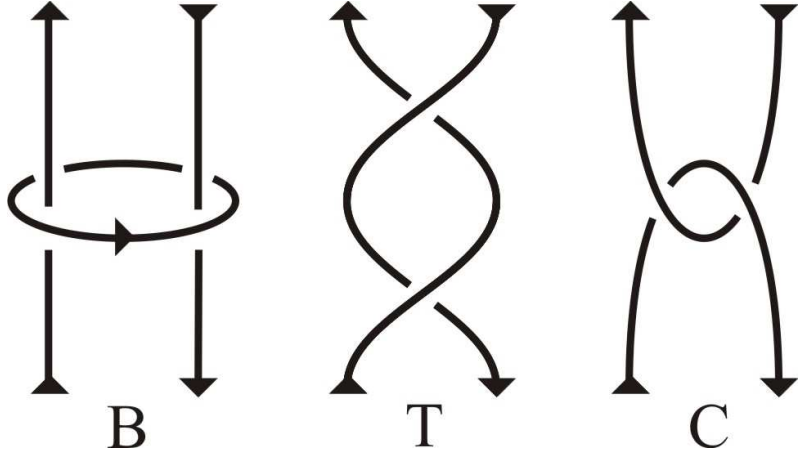}
\caption{The tangles Belt, Twist and Clasp}
\end{center}
\end{figure}

\begin{defs}
Let $a,b,c,d$ be integers such that $b\geq 1$ and $c,d\geq 0$.
Define the Whitehead chain $W_{a,b,c,d}$ to be the closure of the
composition of $a$ tangles of type Twist (or $a$ mirror images of
type Twist when $a<0$), $b$ tangles of type Belt, $c$ tangles of
type Clasp and $d$ mirror images of tangles of type Clasp.
\end{defs}

\noindent The tangles Belt, Twist, Clasp and their mirror images commute, so
the order of composition is immaterial. Therefore the Whitehead
chains are well defined. In the notation of the previous definition the Whitehead link is $W_{0,1,1,0}$ and
the Borromean rings are $W_{0,1,1,1}$.

Our main theorem is the following asymptotic expansion for the
colored Jones polynomial of a Whitehead chain. As usual $x \sim y$ means that the quotient converges to $1$.\\

\noindent \textbf{Main theorem.} \emph{Let $a,b,c,d$ be integers, such that $b\geq 1$ and $c,d\geq 0$.}
$$J_N(W_{a,b,c,d})(e^{\frac{2\pi i}{N}}) \sim \exp\left(\frac{1}{2\pi}\left\{(\mathrm{Vol}(W_{a,b,c,d})+ i \mathrm{CS})N+\mathrm{D}\log(N)+ \mathrm{E}\right\}\right)$$
\begin{enumerate}
\item[-] \emph{The value of $\mathrm{CS}$ is $\frac{-4a+c-d}{8}2\pi^2$ if $c+d = 1$ and $\frac{-4a-7c+7d}{8}2\pi^2$ otherwise.
\item[-] If $b = 1$ we have $\mathrm{D} = 3\pi$ and $\mathrm{E}$ can be expressed explicitly as an integral, see the proof of lemma 4.
\item[-] If $b \geq 2$ the expansion is only valid for odd $N$ and $J_{2M}(W_{a,b,c,d})(e^{\frac{2\pi i}{2M}})=0$. For odd $N$ we find $\mathrm{D} = 2\pi b$ and} $\mathrm{E} =
2\pi(c+d)\log2
+\frac{4a+3c-3d}{4}2\pi^2 i$.
\end{enumerate}

\noindent By taking absolute values and restricting ourselves to the leading term we see that the volume conjecture holds true for $b = 1$, while it is false for $b \geq 2$. In the latter case the volume conjecture is true when we restrict to odd values of $N$. This phenomenon might have something to do with the fact that the complement of the Whitehead chain is hyperbolic when $b = 1$, while the torus decomposition contains a Seifert-fibered piece when $b \geq 2$.

According to the complexified volume conjecture proposed in \cite{MMOTY} the value of $\mathrm{CS}$ is equal to $2\pi^2$ times the
Chern-Simons invariant. For $a = 0$ this is indeed the case because the Chern-Simons invariant is additive with respect to belted sum and its value on the Whitehead link and its mirror image is $\pm\frac{1}{8}$, \cite{Ouyang}.

In \cite{Hikami} the number $\mathrm{D}$ in the above asymptotic
expansion is conjectured to be the number of prime factors of a
knot. In \cite{GM} a different interpretation is given in terms of
the knot complement. In the same paper it is also conjectured that
the number $\mathrm{E}$ is determined by the Ray-Singer torsion of
the complement twisted by the holonomy representation. We hope to
investigate these conjectures for the Whitehead chains in a
subsequent publication.

The main theorem shows that the original volume conjecture may fail even
for non-splittable links, but adds credibility to the following
weaker version of the volume conjecture for links:

\begin{conj}
For any non-splittable link $L$ we have: $$\limsup_{N\to\infty}\frac{2\pi}{N}\log|J_N(L)(e^{\frac{2\pi i}{N}})| = \mathrm{Vol}(L)$$
\end{conj}

\section{Proof of the main theorem} In this section we give an
overview of the proof of the main theorem, postponing the proofs of
the more technical lemmas to the next section.

The first step is to obtain an expression for the colored Jones
polynomials of a general Whitehead chain.

\begin{lem} Let $a,b,c,d$ be integers, with $b \geq 2$ and $c,d \geq 0$. We have the following formulas for the colored Jones polynomial:
$$J_N(W_{a,1,c,d})(e^{\frac{2\pi i}{N}}) =  \phi_N\sum_{n=0}^{N-1}(2n+1)\chi_{N,n}^{4a-c+d}\left(\sum_{k=0}^{N-1-n} S_{n,k}\right)^{c+d}$$
$J_{N}(W_{a,b,c,d})(e^{\frac{2\pi i}{N}})=0$ when $N$ is even and for $N=2M+1$ we have:
$$J_{N}(W_{a,b,c,d})(e^{\frac{2\pi i}{N}}) =  \phi_N N^{b}\chi_{N,M}^{4a-c+d}\left(\sum_{k=0}^{M} S_{M,k}\right)^{c+d}$$
\noindent Where we define $\phi_N =
\exp\left(\frac{(N-1)(c-d)}{N}\pi i\right)$ if $c+d = 1$ and
$\phi_N = (-1)^{(N-1)(c-d)}$ otherwise, and we define $\chi$ and
$S$ by
$$\chi_{N,n}=\exp\left(\frac{n(n+1-N)}{2N}\pi i\right), \quad
S_{n,k} =
\prod_{j=1}^n\frac{2\sin^2(\frac{(k+j)\pi}{N})}{\sin(\frac{j\pi}{N})}$$
\end{lem}

\begin{proof}
Recall that the unnormalized $N$-colored Jones invariants are
intertwining operators of $V_N\otimes V_N$, where $V_N$
corresponds to the $N$-dimensional irreducible
representation of $sl_2$. Using the decomposition $V_N\otimes V_N
= \bigoplus_{n=0}^{N-1}V_{2n+1}$ one can write the unnormalized
colored Jones invariants of the tangles in figure 1 in the
following way \cite{HaoZheng}:
$$\tilde{J}_N(t)=\bigoplus_{n=0}^{N-1}
\mathrm{Tangle}(n,t)\cdot\mathrm{id}_{V_{2n+1}}$$ were the
function $\mathrm{Tangle}(n,t)$ depends on the tangle and
$t=e^{h}$. If we define $[n] =
\frac{t^{\frac{n}{2}}-t^{-\frac{n}{2}}}{t^{\frac{1}{2}}-t^{-\frac{1}{2}}}$
then the functions for the tangles in figure 1 are:

$$\mathrm{T}(n,t) =  t^{n(n+1)},\quad \mathrm{B}(n,t) =  \frac{[N(2n+1)]}{[2n+1]}$$

$$\mathrm{C}(n,t) =  t^{(N^2-1)/2+N(N-1)/2}\sum_{k=0}^{N-1-n}t^{-N(n+k)}\prod_{j=1}^n\frac{(1-t^{N-j-k})(1-t^{j+k})}{1-t^{j}}$$

\noindent The colored Jones polynomials of the mirror-images of these tangles are obtained by replacing $t$ by $t^{-1}$.
The factor $t^{(N^2-1)/2}$ in the formula for $\mathrm{C}(n,t)$ is a correction due to framing that should be included only if both strands in the clasp belong to the same component. For Whitehead chains this means that it should be included only when $c+d = 1$.

By the multiplicativity of the colored Jones invariant with
respect to composition of tangles we can now calculate the colored
Jones polynomial of all Whitehead chains. The general formula for
the normalized version of the colored Jones polynomial of the
Whitehead chains is:
$$J_N(W_{a,b,c,d})(t) = \sum_{n=0}^{N-1}\frac{[2n+1]}{[N]}\mathrm{T}(n,t)^a\mathrm{B}(n,t)^b\mathrm{C}(n,t)^c\mathrm{C}(n,t^{-1})^d$$
The factor $\frac{[2n+1]}{[N]}$ comes from the normalization and from
taking the closure of composition of the tangles.

Routine calculations now yield the above formulas.
\end{proof}

\noindent From now on we will use the shorthand $J_{N,a,b,c,d} =
J_N(W_{a,b,c,d})(e^{\frac{2\pi i}{N}})$.

The next step in proving the main theorem is to investigate the
asymptotics of the above formulas for the colored Jones polynomials
as $N\to \infty$. The factors $S_{n,k}$ turn out to play a crucial
role because they dominate the absolute value of the $(n,k)$-th
term. There exists a unique maximum $\tilde{S}_N =
S_{\left\lfloor\frac{N}{2}\right\rfloor,\left\lfloor\frac{N}{4}\right\rfloor}$
for the $S_{n,k}$. This term dominates all others
so that most of the asymptotics of $J_{N,a,b,c,d}$ can be read off
from this term only. In order to make this precise we need to
compare the other values of $S_{n,k}$ to the maximum value.

Let us choose a fixed number $\delta \in (\frac{1}{2},\frac{c+d+5}{2(c+d+4)})$
once and for all, where the numbers $c$ and $d$ are the parameters
of $W_{a,b,c,d}$. Define $n'=|n-\left\lfloor\frac{N}{2}\right\rfloor|$ and $k' =
|k-\left\lfloor\frac{N}{4}\right\rfloor|$ to be the distances from
the maximum. The next lemma shows how we can estimate the other
values of $S_{n,k}$.

\begin{lem} Using the above definitions of $n',k'$ and
$\delta$ we have:\newline $\mathrm{a)}$ If $n'+k' < N^\delta$ then
$S_{n,k}\tilde{S}_N^{-1} = \exp(-\frac{\pi}{N}(n'^2+2n'k'+2k'^2))
+ \mathcal{O}(N^{3\delta - 2})$, as $N \to \infty$. \newline $\mathrm{b)}$ There are $C,\epsilon >
0$ such that if $n'+k' \geq N^\delta$ then
$S_{n,k}\tilde{S}_N^{-1} < C \exp(-\epsilon N^{2\delta - 1})$
\end{lem}
\noindent In lemma 3 the asymptotics of the maximum value $\tilde{S}_N$ are expressed using the Lobachevski function $\Lambda(\theta) = -\int_0^{\theta}\log|2\sin(x)|dx$.
\begin{lem}
$\tilde{S}_N \sim \exp\left(
\frac{4N}{\pi}\Lambda(\frac{\pi}{4})-\frac{1}{2}\log2N \right), \quad N\to\infty$
\end{lem}
\noindent The asymptotics of $J_{N,a,1,c,d}$ is reduced to that of
$\tilde{S}_N$ by the following lemma. By CS we mean the constant defined in the statement of the main theorem.
\begin{lem} There is a nonzero constant $C \in \mathbb{C}$ such that:
$$J_{N,a,1,c,d}\sim CN^{(c+d+3)/2}\tilde{S}_N^{c+d}e^{\frac{i\mathrm{CS}}{2\pi}},\quad N\to\infty$$
\end{lem}

\noindent There is a similar lemma for the case $b \geq 2$. We confine ourselves to the odd-colored Jones
polynomials, since the even ones are $0$ in $e^{\frac{2\pi
i}{N}}$. In lemma 5 we reduce the asymptotics of the odd colored
Jones polynomial to those of the maximal term $\tilde{S}_N$:

\begin{lem} For $b\geq 2$ and $N$ odd we have: $J_{N,a,b,c,d}\sim $
$$\exp\left(-(c+d)\log\sqrt{2}+(4a+3c-3d)\pi i/4\right) N^{(c+d+2b)/2}\tilde{S}_N^{c+d}e^{\frac{i\mathrm{CS}}{2\pi}}$$
\end{lem}

\noindent Postponing the proofs of these lemmas to the next subsection we
can now prove the main theorem.

\begin{proof} (of the main theorem)\newline
Using an explicit decomposition of the complement into ideal
octahedra it can be shown that $\mathrm{Vol}(W_{a,b,c,d}) =
8(c+d)\Lambda(\frac{\pi}{4})$, see \cite{van der Veen}. Let us
first suppose that $b = 1$. According to lemma 4 there is a
constant $C'$ such that we have: $$J_{N,a,1,c,d}\sim
\exp\left(\log(\tilde{S}_N^{c+d})+\frac{c+d+3}{2}\log(N)+Ni\frac{\mathrm{CS}}{2\pi}+
C' \right)$$ Using lemma 3 we get:
$$\sim \exp\frac{1}{2\pi}\left((\mathrm{Vol}(W_{a,1,c,d})+i\mathrm{CS})N + 3\pi\log(N)+ E\right),\  N\to\infty$$
The case $b \geq 2$ follows in the same way by combining lemma 3 and lemma 5.
\end{proof}

\section{Proof of the lemmas}

In this section we prove the more technical lemmas 2,3,4 and 5.
\begin{proof} (of lemma 2)
The proof of this lemma hinges on the following key estimate of
$S_{n,k}$ in terms of the Lobachevski function $\Lambda(x)$ that was proved
in \cite{HaoZheng}. Define $f(x,y)=
-2\Lambda(x+y)+2\Lambda(y)+\Lambda(x)$. For integers $0\leq n,k, n+k < N$
we have the uniform estimate
$$\log S_{n,k}=\frac{N}{\pi}f(\frac{n\pi}{N},\frac{k\pi}{N})+\mathcal{O}(\log N),\quad N\to \infty$$
Before we can apply this result we first need to show that inside the
triangle $0<x,y,x+y<\pi$ the function $f$ has a unique critical
point $(\frac{\pi}{2},\frac{\pi}{4})$ and reaches its maximum
there, which equals $4\Lambda(\frac{\pi}{4})$. Moreover the Taylor
expansion of $f$ around the critical point is:
$$f(\frac{\pi}{2}+x,\frac{\pi}{4}+y)=f(\frac{\pi}{2},\frac{\pi}{4})-(x^2+2xy+2y^2)+\mathcal{O}(|x|^3+|y|^3)$$
Lemma 2 part a) and b) are direct consequences of these facts once
we note that the difference between $S_{\lfloor
\frac{N}{2}\rfloor,\lfloor \frac{N}{4}\rfloor}$ and the actual
critical value $f(\frac{\pi}{2},\frac{\pi}{4})$ becomes negligibly
small as $N$ grows.

To find the critical points of $f$ in $0<x,y,x+y<\pi$ we use the fundamental theorem of calculus to differentiate $\Lambda$ and
find the system of equations: $2\sin^2(x+y)=\sin(x)$ and
$\sin(x+y)=\sin(y)$. For $0<x,y,x+y<\pi$ this has the unique solution $(x,y)=(\frac{\pi}{2},\frac{\pi}{4})$. To determine the nature of the critical
point we differentiate again and this will be left to the reader.

The value of $f$ at its critical point is
$f(\frac{\pi}{2},\frac{\pi}{4})=-2\Lambda(\frac{3}{4}\pi)+2\Lambda(\frac{\pi}{4})+\Lambda(\frac{\pi}{2})=4\Lambda(\frac{\pi}{4})$
because $\Lambda(x)=-\Lambda(\pi-x)$.
\end{proof}

\noindent The next lemma is an expanded version of a result proven in \cite{HaoZheng}.

\begin{proof} (of lemma 3)
If we define $s_n = -\sum_{j=1}^n \log|2\sin(j\pi/N)|$ then we can write $\log \tilde{S}_N = -2s_{\left\lfloor N/4 \right\rfloor + \left\lfloor N/2 \right\rfloor}+2s_{\left\lfloor N/4 \right\rfloor}+s_{\left\lfloor N/2 \right\rfloor}$. It was shown in \cite{HaoZheng} that for $0<n<\frac{5}{6}N$ we have $s_n = \frac{N}{\pi}\Lambda(\frac{n\pi}{N})-\frac{1}{2}\log n + \mathcal{O}(1)$, as $N\to \infty$. To prove the lemma we need to expand a little further. Assuming that $r = \lim_{N\to\infty}n/N$ exists, we show that $$s_n = \frac{N}{\pi}\Lambda(\frac{n\pi}{N})-\frac{1}{2}\log \frac{n\sin(r\pi)}{r\pi} + \mathcal{O}(N^{-1}), \quad N\to \infty$$
By applying the above expansion for $s_n$ three times we then find the desired expansion for $\tilde{S}_N$.

To prove the expansion for $s_n$ we reason as follows: for $x\in [0,\frac{5}{6}\pi]-\{0,u\}$ we have $$\frac{\sin(x-u)}{x-u} \frac{x}{\sin x} = 1+\frac{-x\cos x +\sin x}{x\sin x}u + \mathcal{O}(u^2), \quad u \to 0$$ Now taking the logarithm and expanding $\log(1+f)$ around $f = 0$ we find: $$\log|\frac{\sin(x-u)}{\sin x}| = \log|\frac{x-u}{x}| +\frac{-x\cos x +\sin x}{x\sin x}u + \mathcal{O}(u^2)$$
Following \cite{HaoZheng} p.7 we set $x = j\pi/N$ and find:
$$s_n - \frac{N}{\pi}\Lambda(\frac{n\pi}{N})= \sum_{j=1}^n \frac{N}{\pi}\int_{0}^{\frac{\pi}{N}}\log\frac{\sin(j\pi/N-u)}{\sin j\pi/N}du = $$
$$-\frac{1}{2}\log n - \frac{1}{2}\log 2\pi + \sum_{j=1}^n \frac{N}{\pi}\int_{0}^{\frac{\pi}{N}}\frac{-(j\pi/N)\cos(j\pi/N) +\sin j\pi/N}{j\pi/N\sin j\pi/N}u du +\mathcal{O}(N^{-1})$$
The term $-\frac{1}{2}\log 2\pi$ is the constant contribution of the Sterling series used in \cite{HaoZheng} p.7. After integration with respect to $u$ we can write the above sum as a Riemann sum. A computation then shows that the limit of this sum equals $\frac{1}{2}(\log r\pi-\log\sin r\pi)$.
\end{proof}

\noindent Although the proof of lemma 4 below is quite long the
main idea is simple: Use lemma 2 to estimate the value of the
colored Jones polynomial in terms of the maximum value
$\tilde{S}_N$.

\begin{proof} (of lemma 4)
For convenience we will assume throughout the proof that $c+d \geq 2$, so that the value of CS is $\frac{-4a-7c+7d}{8}2\pi^2$. The proof in the case $c+d =1$ is completely analogous. According to lemma 1 the formula for the N-colored Jones polynomial of the one-belted
Whitehead chain is
$$J_{N,a,1,c,d} =  \phi_N\sum_{n=0}^{N-1}(2n+1)\chi_{N,n}^{4a-c+d}\left(\sum_{k=0}^{N-1-n} S_{n,k}\right)^{c+d}$$
Define the quotient $Q_N =
J_{N,a,1,c,d}N^{-(c+d+3)/2}\tilde{S}_N^{-(c+d)}e^{\frac{4a+7c-7d}{8}N\pi
i}$. We aim to show that $Q_\infty = \lim_{N\to\infty}Q_N = C$ for
some nonzero complex constant.\newline

\noindent\textbf{Step 1:} By expanding the $c+d$-th power of the
sum we obtain a multi-sum over all N-tuples of natural numbers
$n,\mathbf{k}<N$ such that $n+k_j<N$. Define Central to be the set
of all such tuples that satisfy $n+ \mathrm{max}\ k_j < N^\delta$
and define Far to be the set of tuples such that $n+ \mathrm{max}\
k_j \geq N^\delta$. We can then rewrite $Q_N$ as follows:
$$Q_N=\phi_N e^{\frac{4a+7c-7d}{8}N\pi i} N^{-(c+d+3)/2}\sum_{\mathrm{Central}}(2n+1)\chi_{N,n}^{4a-c+d}S_{n,k_1}\cdots S_{n,k_{c+d}}\tilde{S}_N^{-(c+d)}$$ $$+\phi_N e^{\frac{4a+7c-7d}{8}N\pi i} N^{-(c+d+3)/2}\sum_{\mathrm{Far}}(2n+1)\chi_{N,n}^{4a-c+d}S_{n,k_1}\cdots S_{n,k_{c+d}}\tilde{S}_N^{-(c+d)}$$
The Far sum converges to zero absolutely, because each tuple
$n,\mathbf{k}$ in the sum contains a $k_m$ such that $n+k_m\geq
N^\delta$. According to lemma 2b this means that
$S_{n,k_m}\tilde{S}_N^{-1}=\mathcal{O}(\exp(-\epsilon
N^{2\delta-1}))$. We can estimate all other factors
$S_{n,k_j}\tilde{S}_N^{-1}$ by a constant and conclude that the
sum of absolute values of the Far sum converges to zero since $\delta > 1/2$. So far we have shown that $Q_\infty=$ $$\lim_{N\to\infty} \phi_N e^{\frac{4a+7c-7d}{8}N\pi i} N^{-\frac{c+d+3}{2}}\sum_{\mathrm{Central}}(2n+1)\chi_{N,n}^{4a-c+d}S_{n,k_1}\cdots S_{n,k_{c+d}}\tilde{S}_N^{-(c+d)}$$

\noindent \textbf{Step 2:} In the next step we replace the factor $(2n+1)$
in the expression for $Q_\infty$ by a factor $N$. This is done by
showing that the difference converges to zero. The absolute value
of the difference is less than
$$N^{-(c+d+3)/2}\sum_{\mathrm{Central}}|2n+1-N|S_{n,k_1}\cdots
S_{n,k_{c+d}}\tilde{S}_N^{-(c+d)}$$ For central tuples
$n,\mathbf{k}$ we have $n'<N^\delta$ and hence
$2n+1=N+\mathcal{O}(N^\delta)$. The number of terms in the sum is
of order $\mathcal{O}(N^{(c+d+1)\delta})$ and the product
$S_{n,k_1}\cdots S_{n,k_{c+d}}\tilde{S}_N^{-(c+d)}$ is $\mathcal{O}(1)$ by lemma 2a. Therefore the absolute value of the
difference is at most of order
$\mathcal{O}(N^{-(c+d+3)/2+(c+d+2)\delta})$. Since $\delta <
\frac{c+d+5}{2(c+d+4)}<\frac{c+d+3}{2(c+d+2)}$ the difference
converges to zero and we have: $$Q_\infty= \lim_{N\to\infty}\phi_N
e^{\frac{4a+7c-7d}{8}N\pi i}
N^{-(c+d+1)/2}\sum_{\mathrm{Central}}\chi_{N,n}^{4a-c+d}S_{n,k_1}\cdots
S_{n,k_{c+d}}\tilde{S}_N^{-(c+d)}$$

\noindent \textbf{Step 3:} Define the $(c+d+1)$-variable Gaussian function
$$g_N(x,\mathbf{y}) = \exp-\frac{\pi}{N}\left(( c+d
-(4a+c-d)i)\frac{x^2}{2}+2\sum_{j=1}^{c+d}(y_j+\frac{x}{2})^2\right)$$
and let $\psi = \exp(\frac{4a+3c-3d}{4}\pi i)$. We will show that
$$Q_\infty  = \lim_{N \to \infty}\psi
N^{-(c+d+1)/2}\sum_{\mathrm{Central}}g_N(n',\mathbf{k'})$$
Starting with the sum from step 2 we use lemma 2 to replace the
factors $S_{n,k_j}\tilde{S}^{-1}$ by exponentials and error terms:
$$\lim_{N\to\infty}\phi_N e^{\frac{4a+7c-7d}{8}N\pi i} N^{-(c+d+1)/2}\sum_{\mathrm{Central}}\chi_{N,n}^{4a-c+d}\cdot$$ $$\prod_{j=1}^{c+d}\left(\exp(-\frac{\pi}{N}(n'^2+2n'k_j' + 2k_j'^2))+\mathcal{O}(N^{3\delta-2})\right)$$ To show that the error terms can be removed we estimate the contribution of their absolute values. There are $\mathcal{O}(N^{(c+d+1)\delta})$ terms in the sum so their contribution is of order $\mathcal{O}(N^{(c+d+1)\delta+3\delta-2-(c+d+1)/2})$. This converges to zero because $\delta <\frac{c+d+5}{2(c+d+4)}$. \newline
Next we look at the phase factors $\phi_N$ and $\chi_{N,n}$. We
have $\lim_{N\to \infty}\phi_N =e^{(c-d)\pi i}$. Furthermore:
$$\chi_{N,n}e^{\frac{N\pi i}{8}} =
\exp(\frac{n'^2}{2N}-\frac{N}{8}+\frac{n}{2N}+\frac{N}{8})\pi i =
\exp(\frac{n'^2}{2N}\pi i) \exp(\frac{n}{2N}\pi i)$$ Since
$n=\frac{N}{2}+\mathcal{O}(N^\delta)$ we have
$\lim_{N\to\infty}\frac{n}{2N} = \frac{1}{4}$.  If we apply this
to the above sum we get the desired expression.
\newline

\noindent\textbf{Step 4:} It seems natural to replace the previous
sum by an integral, this is done in step 4. We will show that
$$Q_\infty = \lim_{N \to \infty}\psi N^{-(c+d+1)/2}\int_{|x|+
\mathrm{max}|y_j|<N^\delta-1}g_N(|x|,\mathbf{|y|})\mathrm{d}x\mathrm{d}\mathbf{y}$$
The strategy is to estimate the absolute value of the difference
between sum and integral. First we write the sum from step 3 as
follows:
$$Q_\infty=\lim_{N\to\infty}\psi N^{-(c+d+1)/2}\sum_{\mathrm{Central}}\int_{B(n,\mathbf{k})}g_N(n',\mathbf{k'})\mathrm{d}x\mathrm{d}\mathbf{y}$$ where $B(n,\mathbf{k})=[n-\frac{N}{2}-\frac{1}{2},n-\frac{N}{2}+\frac{1}{2})\times\prod_{j = 1}^{c+d}[k_j-\frac{N}{4}-\frac{1}{2},k_j-\frac{N}{4}+\frac{1}{2})$. The absolute value of the difference between the sum and the proposed integral above is: $$N^{-(c+d+1)/2}|\sum_{\mathrm{central}}\int_{B(n,\mathbf{k})}g_N(n',\mathbf{k'})\mathrm{d}x\mathrm{d}\mathbf{y}-\int_{|x|+ \mathrm{max}|y_j|<N^\delta-1}g_N(|x|,\mathbf{|y|})\mathrm{d}x\mathrm{d}\mathbf{y}|$$
The union of the disjoint blocks $B(n,\mathbf{k})$ as
$n,\mathbf{k}$ runs through Central, covers the entire integration
domain $|x|+ \mathrm{max}|y_j|<N^\delta-1$ so we can subtract the
integrals. However some blocks continue over the boundary of the
domain of integration, resulting in a slight error. The terms
$n,\mathbf{k}$ such that the corresponding blocks
$B(n,\mathbf{k})$ intersect the complement of the domain will be
called Border terms. We can estimate the above quantity as
follows: $$<
N^{-(c+d+1)/2}\sum_{\mathrm{Central}}\int_{B(n,\mathbf{k})}|g_N(n',\mathbf{k'})-g_N(|x|,\mathbf{|y|})|\mathrm{d}x\mathrm{d}\mathbf{y}$$
$$+N^{-(c+d+1)/2}\sum_{\mathrm{Border}}|g_N(n',\mathbf{k'})|$$
Both sums will be shown to converge to zero, we start with the
second one. Since the number of terms on the Border of the
integration domain is $\mathcal{O}(N^{(c+d)\delta})$ and
$g_N(n',\mathbf{k'})= \mathcal{O}(1)$ the second sum is of order
$\mathcal{O}(N^{-(c+d+1)/2+(c+d)\delta})$. This implies that the
sum converges to zero because
$\delta<\frac{c+d+5}{2(c+d+4)}<\frac{c+d+1}{2(c+d)}$.\newline For
the first sum we need to estimate the integrands:
$$|g_N(n',\mathbf{k'})-g_N(|x|,\mathbf{|y|})|=|g_N(n',\mathbf{k'})||1-$$ $$\exp\frac{-\pi}{N}\left((c+d -(4a+c-d)i)\frac{|x|^2-n'^2}{2}+2\sum_{j=1}^{c+d}(|y_j|+\frac{|x|}{2})^2-(k_j'+\frac{n'}{2})^2\right)|$$ The last expression is of order $\mathcal{O}(N^{\delta-1})$ because $|g_N(n',\mathbf{k'})|=\mathcal{O}(1)$ and $|x|^2-n'^2 = (|x|+n')(|x|-n')$ and $(|y_j|+\frac{|x|}{2})^2-(k_j'+\frac{n'}{2})^2 = (|y_j|+\frac{|x|}{2}+k_j'+\frac{n'}{2})(|y_j|+\frac{|x|}{2}-k_j'-\frac{n'}{2})$. We integrate over $B(n,\mathbf{k})$ so $|n'-|x||<1$ and $|k_j'-|y_j||<1$ and all terms are central so $|x|+n'<N^\delta$ and $k_j'+|y_j|<N^\delta$. Therefore the second sum is of order $\mathcal{O}(N^{-(c+d+1)/2+\delta-1+(c+d+1)\delta})$ and thus converges to zero, because $\delta<\frac{c+d+5}{2(c+d+4)}<\frac{c+d+3}{2(c+d+2)}$.
\newline

\noindent\textbf{Step 5:} We make the substitution $x=\sqrt{N}w$
and $y_j = \sqrt{N}z_j$ with Jacobian $N^{(c+d+1)/2}$. This gives:
$$Q_\infty = \lim_{M \to \infty}\psi
\int_{w+\mathrm{max}|z_j|<\frac{(N^\delta-1)}{\sqrt{N}}}$$
$$\exp(-\pi\left(( c+d
-(4a+c-d)i)\frac{|w|^2}{2}+2\sum_{j=1}^{c+d}(|z_j|+\frac{|x|}{2})^2\right))\mathrm{d}w\mathrm{d}\mathbf{z}$$
Now $\delta>\frac{1}{2}$ and the integrand is rapidly decreasing
so the limit exists and is equal to:
$$\psi \int_{\mathbb{R}^{c+d+1}}\exp(-\pi\left(( c+d -(4a+c-d)i)\frac{|w|^2}{2}+2\sum_{j=1}^{c+d}(|z_j|+\frac{|w|}{2})^2\right))\mathrm{d}w\mathrm{d}\mathbf{z}$$

\noindent\textbf{Step 6:} In this final step we need to show that
$Q_\infty \neq 0$. We know that $\psi \neq 0$, but what about the
complicated Gaussian integral above? We can write it as an
iterated integral and get rid of the absolute value signs by
integrating $2^{c+d+1}$ times over the positive hyper-quadrant
$w,z_j>0$: $$Q_\infty= \psi 2^{c+d+1}\int_{0}^\infty\exp(-\pi( c+d
-(4a+c-d)i)\frac{w^2}{2})\mathrm{d}w\ \cdot$$
$$\prod_{j=1}^{c+d}\int_{0}^\infty\exp(-2\pi(z_j+\frac{w}{2})^2))\mathrm{d}z_j\mathrm{d}w$$
Using the substitutions $y_j = \sqrt{2\pi}(z_j+\frac{w}{2})$ we
get the integral below. We only want to check that the integral is
nonzero so all constants in front of the integral and the Jacobian
are written as $C$. $$Q_\infty= C\int_{0}^\infty\exp(-\pi( c+d
-(4a+c-d)i)\frac{w^2}{2})\cdot\prod_{j=1}^{c+d}\int_{\frac{\sqrt{\pi}w}{\sqrt{2}}}^\infty\exp(-y_j^2)\mathrm{d}y_j\mathrm{d}w$$
If we define the complementary error function by $\mathrm{erfc}(x)
= \frac{2}{\sqrt{\pi}}\int_{x}^\infty \exp(-y^2)\mathrm{d}y$, then
we can write the integral more concisely as:
$$Q_\infty= C\int_{0}^\infty\exp(-\pi( c+d -(4a+c-d)i)\frac{w^2}{2})\mathrm{erfc}^{c+d}\left(\frac{\sqrt{\pi}w}{\sqrt{2}}\right)\mathrm{d}w =$$
$$C\int_{0}^\infty\exp(-\pi(c+d)\frac{w^2}{2})\mathrm{erfc}^{c+d}\left(\frac{\sqrt{\pi}w}{\sqrt{2}}\right)\cdot$$ $$\left(\cos(\frac{4a+c-d}{2}\pi w^2) + i\sin(\frac{4a+c-d}{2}\pi w^2)\right)\mathrm{d}w $$
After the substitution $x = \sqrt{\frac{4a+c-d}{2}}w$ we note that
the imaginary part of integral is nonzero, because both the
exponential function and the complementary error function are
bounded and decrease monotonically and the decrease is quite
rapid. Therefore the integral over the part where the sine is
positive dominates the integral over its complement, showing that
the integral is a positive scalar multiple of $C$.
\end{proof}

The same proof also works for lemma 5 but in this case it is much
easier because the value of $n$ is fixed at the maximum. The Gaussian integral at the end of step 5 is now a standard Gaussian integral so that the constant in front of lemma 5 can be computed easily.

\end{document}